\pdfoutput=1
\RequirePackage{ifpdf}
\ifpdf 
\documentclass[pdftex]{sigma}
\else
\documentclass{sigma}
\fi

\usepackage{multirow}
\usepackage{array}
\newcolumntype{I}{!{\vrule width 2pt}}
\newlength\savedarrayrulewidth
\newcommand\boldhline{\noalign{\global\savedarrayrulewidth\arrayrulewidth
  \global\arrayrulewidth 2pt}%
  \hline
  \noalign{\global\arrayrulewidth\savedarrayrulewidth}}

\DeclareMathOperator{\SLtwo}{SL_2}

\begin{document}


\renewcommand{\thefootnote}{$\star$}

\renewcommand{\PaperNumber}{075}

\FirstPageHeading

\ShortArticleName{Sylvester versus Gundelf\/inger}

\ArticleName{Sylvester versus Gundelf\/inger\footnote{This
paper is a contribution to the Special Issue ``Symmetries of Dif\/ferential Equations: Frames, Invariants and Applications''. The full collection is available at \href{http://www.emis.de/journals/SIGMA/SDE2012.html}{http://www.emis.de/journals/SIGMA/SDE2012.html}}}

\Author{Andries E.~BROUWER~$^\dag$ and Mihaela POPOVICIU~$^\ddag$}

\AuthorNameForHeading{A.E.~Brouwer and M.~Popoviciu}

\Address{$^\dag$~Department of Mathematics and Computer Science, Technische Universiteit Eindhoven,\\
\hphantom{$^\dag$}~P.O.~Box 513, 5600 MB Eindhoven, The Netherlands}
\EmailD{\href{mailto:aeb@cwi.nl}{aeb@cwi.nl}}
\URLaddressD{\url{http://www.win.tue.nl/~aeb/}}

\Address{$^\ddag$~Mathematisches Institut, Universit\"at Basel, Rheinsprung 21,
CH-4051 Basel, Switzerland}
\EmailD{\href{mailto:mihaela.popoviciu@unibas.ch}{mihaela.popoviciu@unibas.ch}}

\ArticleDates{Received July 18, 2012, in f\/inal form October 12, 2012; Published online October 19, 2012}

\Abstract{Let $V_n$ be the $\SLtwo$-module of binary forms of degree $n$
and let $V = V_1 \oplus V_3 \oplus V_4$.
We show that the minimum number of generators of the algebra
$R = \mathbb{C}[V]^{\SLtwo}$ of polynomial functions
on $V$ invariant under the action of $\SLtwo$
equals 63. This settles a 143-year old question.}


\Classification{13A15; 68W30}

\vspace{-2mm}

\renewcommand{\thefootnote}{\arabic{footnote}}
\setcounter{footnote}{0}

\section{Introduction}

In 1868 Gordan~\cite{Go} proved that the algebra of invariants of
binary forms of given degree is f\/initely generated.
This came as a surprise to Cayley and Sylvester, who had believed
that the quintic and septimic had covariant resp. invariant rings
that were not f\/initely generated.
\begin{quote}
{\it The number of invariants is first infinite in the case of a quantic
of the seventh order, or septic; the number of covariants is first infinite
in the case of a quantic of the fifth order, or quintic.}
(Cayley \cite{Cay})
\end{quote}

However, f\/inding a minimal set of generators for these algebras is even today
an open problem in all but a few small cases.
In the case of $V_4 \oplus V_4$, Gordan found a generating set
of size 30, and Sylvester \cite{Sylv0} showed that two of these generators
are superf\/luous. He added
\begin{quote}
{\it J'ajouterai seulement que cette preuve \'eclatante de l'insuffisance
de la m\'ethode de M. Gordan et de son \'ecole, pour s\'eparer les
formes v\'eritablement \'el\'ementaires des formes superflues qui
s'y rattachent $($insuffisance reconnue par M. Gordan lui-m\^eme
de la mani\`ere la plus loyale dans son discours inaugural
prononc\'e \`a Erlangen$)$, n'\^ote rien \`a la valeur immense
du service qu'il a rendu \`a l'Alg\`ebre, en ayant le premier
d\'emontr\'e l'existence d'une limite au nombre de ces formes.}
\end{quote}

This note focuses on the covariants of $V_3\oplus V_4$, a
case which illustrates the controversy between the German
and English schools in the 19th century. The German school,
following Clebsch and Gordan, was able to construct a system
of generators for the algebra of invariants of binary forms,
with no guarantee that the system was minimal. The English school, following
Cayley and Sylvester, aimed to determine the number of independent
generators. Sylvester used in his computations his `fundamental postulate'
(not def\/ined here),
which turned out to hold only in small cases.
Counterexamples were given by Hammond~\cite{Ha} and Morley~\cite{Mo}.

In 1869 Gundelf\/inger, a student of Clebsch, wrote a thesis~\cite{Gund}
where he constructed generators for the covariants of $V_3 \oplus V_4$
`in ordinary symbolic notation', after Clebsch had given him this system
as computed by Gordan in his `obscure' notation
(cf.~\cite[pp.\,270--272]{GordanS}).
He found 20 generators for the invariants and 64 for the
covariants.

Sylvester used the Poincar\'e series together with his
fundamental postulate to show that there could be only 61
independent generators for the covariants of $V_3 \oplus V_4$,
and wrote a series of papers~\cite{Sylv3,Sylv2,Sylv1,SylvF}
showing the superiority of the English methods over the German.

In the f\/irst paper Sylvester uses his method (which he calls `{\em tamisage}')
to derive the numbers of generators of given degrees in the coef\/f\/icients
of~$V_4$ and those of~$V_3$, and given order in the variables~$x$,~$y$.
The following table is taken from \cite{SylvF}:
\begin{figure}[h!]\centering
\includegraphics[scale=0.56]{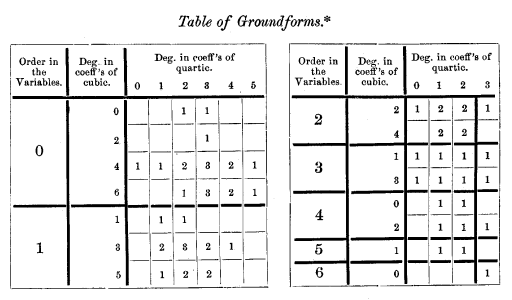}\\
\vspace{1mm}
\includegraphics[scale=0.6]{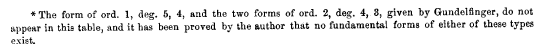}
\end{figure}

In the second paper he observes that it follows from the Poincar\'e series
that there are 8 li\-near\-ly independent covariants of order~2 and
multidegree~(4,3). Next, he constructs 8 reducible such covariants
(products of covariants of lower degree) and argues that these are
linearly independent. However, the forms are dependent and only seven
are independent.
He f\/inishes with the announcement
\begin{figure}[h!]\centering
\includegraphics[scale=0.6]{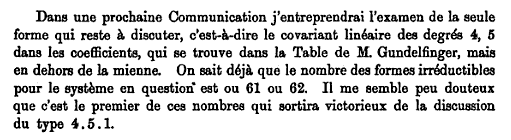}
\end{figure}

In the third paper he observes that it follows from the Poincar\'e series
that there are 12 linearly independent covariants of order 1 and
multidegree (5,4). Next, he constructs 12 reducible such covariants
and argues that these are linearly independent. However, the forms are
dependent and only eleven are independent.
He concludes
\begin{center}
\includegraphics[scale=0.6]{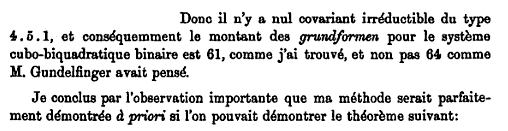}
\\
(false theorem omitted)\\
\includegraphics[scale=0.6]{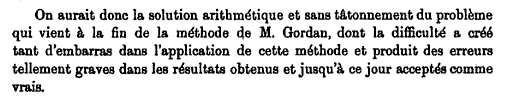}
\end{center}

\noindent
Here we show that the actual minimal number of generators for the
covariants of~$V_3\oplus V_4$ is~63. Our results coincide with those
of Sylvester and Gundelf\/inger, with two exceptions: we show that one
needs one generating covariant of order~1 and multidegree~(5,4),
and (only) one generating covariant of order~2 and multidegree~(4,3).

For completeness we give the corrected version of Sylvester's table.
The two corrected entries are underlined.

\smallskip\begin{center}{\footnotesize
\begin{tabular}{|c|c|cccccc|}
\hline
        &  {Deg. in}   & \multicolumn{6}{c|}{Deg. in coef\/f's}\\
{Order} & {coef\/f's of} & \multicolumn{6}{c|}{of quartic.}  \\
        &  {cubic.}    & 0 & 1 & 2 & 3 & {4} & {5} \\
\boldhline
\multirow{4}{*}{0} & 0 &   &   & 1 & 1 &   &   \\
                   & 2 &   &   &   & 1 &   &   \\
                   & 4 & 1 & 1 & 2 & 3 & 2 & 1 \\
                   & 6 &   &   & 1 & 3 & 2 & 1 \\
\boldhline
\multirow{3}{*}{1} & 1 &   & 1 & 1 &   &   & \\
                   & 3 &   & 2 & 3 & 2 & 1 & \\
                   & 5 &   & 1 & 2 & 2 & \underline{1} & \\
\hline
\end{tabular}
\quad\begin{tabular}{|c|c|cccc|}
\hline
      &   Deg. in  & \multicolumn{4}{c|}{Deg. in coef\/f's}\\
Order & coef\/f's of & \multicolumn{4}{c|}{of quartic.}\\
      &   cubic.   & 0 & 1 &  2 & 3 \\
\boldhline
\multirow{2}{*}{2} & 2 & 1 & 2 & 2 & 1 \\
                   & 4 &   & 2 & 2 & \underline{1} \\
\boldhline
\multirow{2}{*}{3} & 1 & 1 & 1 & 1 & 1 \\
                   & 3 & 1 & 1 & 1 & 1 \\
\boldhline
\multirow{2}{*}{4} & 0 &   & 1 & 1 &   \\
                   & 2 &   & 1 & 1 & 1 \\
\boldhline
                 5 & 1 &   & 1 & 1 &   \\
\boldhline
                 6 & 0 &   &   &   & 1 \\
\hline
\end{tabular}
}\end{center}

\subsection{Invariants and covariants}
Let $V$ be a f\/inite-dimensional vector space over a f\/ield~$k$, with basis
$e_1,\ldots,e_m$. Let $x_i$ be the coordinate function
def\/ined by $x_i(\sum \xi_h e_h) = \xi_i$.
The algebra $k[V]$ of polynomial functions on~$V$ is by def\/inition
the algebra generated by the~$x_i$. (It does not depend on the choice
of basis $e_1,\ldots,e_m$.)
Let~$G$ be a group of linear transformations of~$V$.
It acts on~$k[V]$ via the action $(g\cdot f)(v) = f(g^{-1}v)$.
Invariant theory studies $k[V]^G$, the algebra of $G$-invariant
polynomial functions on $V$, i.e., the $f \in k[V]$ such that
$g\cdot f = f$ for all $g \in G$.

A covariant of order $m$ and degree $d$ of $V$ is a $G$-equivariant
homogeneous polynomial map $\phi : V\rightarrow V_m$ of degree $d$. In
other words, $\phi (g \cdot v)=g \cdot \phi (v)$, for all $g \in G$,
and $\phi (tv)=t^d \phi (v)$, for all $t\in k$. In particular,
the covariants of $V$ of order~0 are the invariants of~$V$.

Below we shall take $k = \mathbb{C}$, $G = \SLtwo(k)$, and
$V = V_{n_1} \oplus \cdots \oplus V_{n_p}$, where
$V_n$ is the vector space (of dimension $n+1$) consisting of 0
and the binary forms of degree $n$, that is,
of the homogeneous polynomials of degree $n$
\[
v(x,y) = a_0 x^n + \binom{n}{1} a_1 x^{n-1} y + \cdots +\binom{n}{n-1}
a_{n-1} x y^{n-1}+ a_n y^n,
\]
in two variables. This $V_n$ is the $n$-th graded part of $k[W]$,
where $W$ is a 2-dimensional vector space over $\mathbb{C}$ with natural action
of $\SLtwo$, hence has a natural action of $\SLtwo$.

The main way to construct covariants is via {\em transvectants}
(\"Uberschiebungen). These are derived from the Clebsch--Gordan
decomposition of the $\SLtwo$-module $V_m\otimes V_n$, with $m\geq n$:
\[
V_m\otimes V_n\simeq V_{m+n}\oplus V_{m+n-2}\oplus \cdots \oplus V_{m-n}.
\]
This decomposition def\/ines for each $p$, $0\leq p \leq n$,
an $\SLtwo $-equivariant linear map $V_m\otimes V_n \to V_{m+n-2p}$,
denoted $(g,h)\mapsto (g,h)_p$, and called the {\it p-th transvectant}.
It is given explicitly by the following formula:
\[
(g,h) \mapsto (g,h)_p:=\frac{(m-p)!(n-p)!}{m!n!}
\sum_{i=0}^p (-1)^i \binom{p}{i}
\frac{\partial ^p g}{\partial x^{p-i}\partial y^i}
\frac{\partial ^p h}{\partial x^i \partial y^{p-i}}
\]
(see \cite[Chapter 5]{Olver}).

The covariants of $V$ can be identif\/ied with the invariants of
$V_1\oplus V$: we have $V_1\oplus V\simeq V_1^*\oplus V$ as
$\SLtwo$-representations and the set of covariants of $V$ is isomorphic to
$k[V_1^*\oplus V]^{\SLtwo}$ (see \cite[Chapter 15]{Pro07}). Each
covariant $\phi $ of $V$ of order $m$ corresponds to
the invariant of $V_1\oplus V$ def\/ined by the transvectant $(\phi (v)
,l ^m)_m$, where $l \in V_1$.

\section[The generators of the invariants of $V_1\oplus V_3 \oplus V_4$]{The generators of the invariants of $\boldsymbol{V_1\oplus V_3 \oplus V_4}$}

We identify the covariants of $V_3\oplus V_4$ with the invariants
of $V_1\oplus V_3 \oplus V_4$ and show that a minimal set of generators
for the algebra of invariants of this module has size 63.

Doing this type of work requires f\/inding dependencies.
Gundelf\/inger did not try to do this exhaustively, but following Gordan
he only noted the obvious ones. Sylvester tried, and made
some mistakes, no doubt because he already knew what answer he wanted.
For us this is relatively easy~-- a~modern computer has no
problems computing the rank of a 40000 by 600000 matrix (which is what
is needed in the most straightforward approach).

We had a dif\/ferent problem: up to which degree should we compute
covariants or inva\-riants? Gundelf\/inger `just' followed Gordan's
algorithm, but as far as we know that has not been implemented yet.

The secret knowledge known today but not in the 19th century,
is that the ring $R$ of invariants of $V_1 \oplus V_3 \oplus V_4$
(or any such ring) is Cohen--Macaulay (see~\cite{HoRo}).
It has a homogeneous system
of parameters (hsop) $j_1,\ldots,j_r$, algebraically independent,
and f\/initely many further generators $i_1,\ldots,i_s$, such that
every invariant can be uniquely written as a linear combination
of products $i_m j_{m_1} \cdots j_{m_h}$.
It follows that the Poincar\'e series $P(t) = \sum d_i t^i$,
where $d_i$ is the dimension of the degree $i$ part of $R$, is of the form
\[
P(t) = \frac{t^{a_1}+\cdots+t^{a_s}}{(1-t^{b_1})\cdots(1-t^{b_r})},
\]
where the $a_h$ and $b_h$ are the degrees of the $i_h$ and $j_h$.

For the module $V = V_{n_1} \oplus \cdots \oplus V_{n_p}$,
with $\sum_i n_i \ge 3$, one has
$P(t^{-1}) = (-1)^{d-3} t^d P(t)$ where $d = \sum_i (n_i+1)$
(by \cite[Corollary~2]{Springer80} for $p=1$,
and by \cite[Theorem~2]{Brion} in general),
so that $\max_h a_h - \sum_j b_j = - \sum_i (n_i+1)$. Therefore,
in order to f\/ind $\max_h \{a_h, b_h\}$ it suf\/f\/ices to f\/ind the~$b_h$.

The power series $P(t) = \sum d_i t^i$ is known from
Cayley \cite{Cay} and Sylvester \cite{Sylv4} (cf.~\cite[3.3.4]{Springer77}).
In the present case,
\begin{gather*}
P(t) =   1+t^2+2t^3+5t^4+10t^5+18t^6+31t^7+55t^8+92t^9 +
  144t^{10}+223t^{11}+341t^{12}\\
\hphantom{P(t) =}{}
+499t^{13}+725t^{14}+1031t^{15} +
  1436t^{16}+1978t^{17}+2685t^{18}+3592t^{19}+4761t^{20}\\
\hphantom{P(t) =}{}
  +6235t^{21} +
  8078t^{22}+10379t^{23}+13226t^{24}+16698t^{25}+20937t^{26} +
  26069t^{27}\\
\hphantom{P(t) =}{}
  +32230t^{28}+39614t^{29}+48401t^{30}+ \cdots   \\
\hphantom{P(t)}{}
=   \frac{a(t)}{(1-t^3)(1-t^4)^2(1-t^5)^2(1-t^6)^2(1-t^7)},
\end{gather*}
where
\begin{gather*}
a(t) =   1+t^2+t^3+3t^4+7t^5+12t^6+21t^7+32t^8+47t^9 +
  58t^{10}+72t^{11}+83t^{12}+89t^{13}\\
\hphantom{a(t) =}{}
  +94t^{14}+94t^{15}+89t^{16} +
  83t^{17}+72t^{18}+58t^{19}+47t^{20}+32t^{21}+21t^{22}+12t^{23}\\
\hphantom{a(t) =}{}
   +
  7t^{24}+ 3t^{25}+t^{26}+t^{27}+t^{29},
\end{gather*}
and it follows that computing invariants up to degree 29 suf\/f\/ices
if we know that there is a hsop with degrees 3, 4, 4, 5, 5, 6, 6, 7.

\subsection{Finding a hsop}

Hilbert introduced in the 19th century the notion of
{\em nullcone}. If $V$ is an $\SLtwo $-module, then the nullcone
$\mathcal{N} (V)$ of $V$ is the set of elements of $V$ on which all invariants
of $V$ of positive degree vanish.
The elements of~$\mathcal{N} (V)$ are called {\em nullforms}.
One can show that a binary form $f\in V_n$
is a nullform if and only if $f$ has a root of multiplicity $>
\frac{n}{2}$ (this is a consequence of the Hilbert--Mumford
criterion, see \cite[\S~2.4.1]{DeKe}). Similarly, if we have
$p$ binary forms $f_1,\ldots, f_p$ of degrees $n_1,\ldots ,n_p$,
then $(f_1,\ldots,f_p)\in \mathcal{N} (V_{n_1}\oplus \cdots \oplus V_{n_p})$
if and only if $f_1,\ldots , f_p$ have a~common root that has
multiplicity $>\frac{n_i}{2}$ in $f_i$, for all $i=1,\ldots ,p$.
In our particular case, if $(l,c,q)\in V_1\oplus V_3 \oplus V_4$, then
$(l,c,q)\in \mathcal{N} (V_1\oplus V_3 \oplus V_4)$ if and only if
$l^2\,|\, c$ and $l^3\,|\, q$.

Let $\mathcal{V}(J)$ stand
for the vanishing locus of $J$. The following result, due to Hilbert, gives a~characterisation of homogeneous systems of parameters of
$k[V_{n_1}\oplus \cdots \oplus V_{n_p}]^{\SLtwo }$ as sets
that def\/ine the nullcone of $\mathcal{N} (V_{n_1}\oplus \cdots \oplus V_{n_p})$:
\begin{proposition}[Hilbert \cite{Hi2}] \label{hilbert}
Let $V=V_{n_1}\oplus \cdots \oplus V_{n_p}$, and $R =
k[V]^{\SLtwo }$,
and $m =n_1+\cdots +n_p+p-3>0$. A set $\{j_1,\ldots, j_m\}$
of homogeneous elements of~$R$ is a system of parameters of~$R$ if and
only if $\mathcal{V}(j_1,\ldots, j_m)=\mathcal{N}(V)$.
\end{proposition}

Let our binary forms $l\in V_1$, $c \in V_3$, $q\in V_4$ be
\begin{gather*}
l =c_0 x+c_1y,  \\
c  =a_0 x^3+3 a_1 x^2 y+3a_2 x y^2 +a_3y^3,  \\
q =b_0x^4+4 b_1 x^3y+6b_2 x^2y^2+4 b_3 xy^3+b_4 y^4,
\end{gather*}
and consider the following invariants:
\begin{alignat*}{3}
& k_2 = (q,q)_4, \qquad &&  k_3 = ((q,q)_2,q)_4, &\\
& k_{4,1} = ((c,c)_2,(c,c)_2)_2, \qquad && k_{4,2} = (lc,lc)_4, & \\
& k_{4,3} = (c,l^3)_3,  \qquad && k_{5,1} = \big((q,(q,q)_2)_1,c^2\big)_6, & \\
& k_{5,2} = \big((q,c^2)_2,c^2\big)_6,   \qquad && k_{5,3} = \big(q,l^4\big)_4, & \\
& k_{6,1} = \big([(c,c)_2]^2,(q,q)_2\big)_4,   \qquad && k_{6,2} = ((lc,lc)_2,lc)_4, & \\
& k_{6,3} = \big((q,q)_2,l^4\big)_4,\qquad && k_7 = \big(c^4,q^3\big)_{12}. &
\end{alignat*}
We prove the following
\begin{proposition}
With the notations above, the invariants
\begin{alignat*}{5}
& j_1 = k_3, \qquad && j_2 = k_{4,1}+k_2^2, \qquad && j_3 = k_{4,2}+k_{4,3}-k_2^2, \qquad && j_4 = k_{5,1}+k_{5,2}, & \\
& j_5 = k_{5,3}, \qquad && j_6 = k_{6,1}+k_{6,2}, \qquad && j_7 = k_{6,3}, \qquad && j_8 = k_7, &
\end{alignat*}
$($of degrees $3$, $4$, $4$, $5$, $5$, $6$, $6$, $7$, respectively$)$
form a system of parameters of $k[V_1 \oplus V_3 \oplus V_4]^{\SLtwo }$.
\end{proposition}
\begin{proof}
We show that $\mathcal{V}(j_1, \dots, j_8) = \mathcal{N}(V_1\oplus V_3
\oplus V_4)$.
Consider three cases.

Case 1: $q=0$.
In this case, the vanishing of $j_1, \dots, j_8$ reduces to
$k_{4,1}=k_{4,2}+k_{4,3}=k_{6,2}=0$, which implies that
$(l,c)\in \mathcal{N}(V_1\oplus V_3)$.
Indeed, if $k_{4,1}=0$, then $c$ is a nullform, and, without loss
of generality, we may suppose that $x^2\,|\, c$, i.e.~$a_2=a_3=0$.
But then $k_{6,2} \sim a_1^3 c_1^3$ (we use~`$\sim $' for
equality up to a nonzero constant).
If $c_1=0$, then $(l,c)\in \mathcal{N}(V_1\oplus V_3)$.
If $c_1\neq 0$, then $a_1=0$ and $k_{4,2}+k_{4,3}\sim a_0c_1^3$.
Hence $a_0=0$, so that $c = 0$ and $(l,c)\in \mathcal{N}(V_1\oplus V_3)$.

Case 2: $l=0$.
In this case, the vanishing of $j_1, \dots, j_8$ reduces to
$k_2 = k_3 = k_{4,1} = k_{5,1}+k_{5,2} = k_{6,1} = k_7 = 0$,
which implies that $(c,q)\in \mathcal{N} (V_3\oplus V_4)$.
Indeed, the vanishing of $k_2$, $k_3$, $k_{4,1}$
implies that $c$ and $q$ are nullforms. If~$c$ or~$q$ vanish identically,
then the statement is clear. Otherwise, if the double zero of~$c$ and the
triple zero of $q$ do not coincide, we may suppose, without loss of
generality, that $x^2 \,|\, c$ and $y^3\,|\,  q$,
i.e. $a_2=a_3=b_0=b_1=b_2=0$. Then
$k_{6,1}\sim a_1^4 b_3^2$. If $a_1=0$, then
$k_{5,1}+k_{5,2}\sim a_0^2 b_3^3$,
and $k_7 \sim a_0^4 b_4^3$, which contradicts the assumption $c, q \neq 0$.
If $b_3=0$, then $k_{5,1}+k_{5,2}\sim a_1^4 b_4^{\phantom{0}}$ and
$k_7 \sim a_0^4 b_4^3$, which again contradicts the assumption $c, q \neq 0$.

Case 3: $q,\,l\neq 0$.
In this case, $j_5 = 0$ implies that $q$ and $l$ have a common root
(up to a~constant, $j_5$ is the resultant of~$q$ and~$l$).
Without loss of generality, we can suppose that the common factor of~$q$ and
$l$ is $x$, i.e., $c_1=b_4=0$ and $c_0 \ne 0$.
Then $j_7 \sim b_3^2 c_0^4$, which implies $b_3=0$.
Then $j_1\sim b_2^3$, which implies $b_2=0$.
Then $a_3$ becomes a factor of $j_8$.
If $a_3=0$, then $j_3\sim a_2^2c_0^2$, which implies $a_2=0$, and then
$(l,c,q) \in \mathcal{N}(V_1 \oplus V_3 \oplus V_4)$. If $a_3\neq 0$,
we may take $a_3 = c_0 = 1$. Now
\[
j_3 \sim 3a_2^2-3a_1-2,
\]
and it follows that $a_1= a_2^2-\frac{2}{3}$. Then
\[
j_6 \sim 27a_2^3-54a_2-27a_0-256b_1^2,
\]
and it follows that
  $a_0=a_2^3-2a_2-\frac{256}{27}b_1^2$. Then
\[
j_4 \sim 36b_0-144a_2b_1-949b_1^3,
\]
and it follows that
  $b_0=4a_2b_1+\frac{949}{36}b_1^3$. Then
\[
j_2 \sim 27-2048 b_1^4,\qquad
j_8  \sim b_1^5\big(33205248-4273351745 b_1^4\big).
\]
But $j_2 = j_8 = 0$ has no solution. This settles Case~3.

By Proposition~\ref{hilbert}, it follows that these eight invariants
form a hsop of the ring of invariants of $V_1\oplus V_3\oplus V_4$.
\end{proof}

\subsection{The degrees of the generators}

The Poincar\'e series of the ring of invariants of $V_1\oplus V_3\oplus
V_4$ tells us which is the maximal degree in which we
have to look for generators, namely~29. For each $i \leq 29 $ we do the
following: multiply invariants of smaller degrees to see what part
of the vector space of invariants of degree~$i$ is known. The
Poincar\'e series tells us how big the dimension of this vector space
is, and if the known invariants do not yet span this vector space,
one constructs in some way further invariants, until they do span.
In the following table~$i$ denotes the degree of the generators,
and $d_i$ the number of generators of degree~$i$ needed:
$$
\begin{tabular}{@{\,}c|@{\hspace{2.0mm}}c
@{\hspace{2.0mm}}c
@{\hspace{2.0mm}}c
@{\hspace{2.0mm}}c
@{\hspace{2.0mm}}c
@{\hspace{2.0mm}}c
@{\hspace{2.0mm}}c
@{\hspace{2.0mm}}c
@{\hspace{2.0mm}}c
@{\hspace{2.0mm}}c
@{\hspace{2.0mm}}c@{\,}} $i$ & $2$ & $3$ &$4$ & $5$ &$6$
  &$7$ & $8$ & $9$ & $10$ & $11$ \\
  \hline $d_i$ & $1$ &$2$ &$4$ &$8$  &$10$ &$13$ & $11$ & $10$ & $3$ & $1$
\end{tabular}
$$
For $12 \le i \le 29$ no further generators are needed, and it follows
that the minimal number of generators is 63.

\subsection*{Acknowledgements}

The second author is partially supported by
the Swiss National Science Foundation.

\pdfbookmark[1]{References}{ref}
\LastPageEnding


\begin{thebibliography}{99}
\footnotesize\itemsep=0pt

\bibitem{Brion}
Brion M., Invariants de plusieurs formes binaires, \textit{Bull. Soc. Math.
  France} \textbf{110} (1982), 429--445.

\bibitem{Cay}
Cayley A., A second memoir upon quantics, \textit{Phil. Trans. Royal Soc.
  London} \textbf{146} (1856), 101--126.

\bibitem{DeKe}
Derksen H., Kemper G., Computational invariant theory, \textit{Encyclopaedia of
  Mathematical Sciences}, Vol.~130, Springer-Verlag, Berlin, 2002.

\bibitem{Go}
Gordan P., Beweis, dass jede Covariante und Invariante einer bin\"aren Form
  eine ganze Funktion mit numerischen Coef\/f\/izienten einer endlichen Anzahl
  solcher Formen ist, \href{http://dx.doi.org/10.1515/crll.1868.69.323}{\textit{J.~Reine Angew. Math.}} \textbf{69} (1868),
  323--354.

\bibitem{GordanS}
Gordan P., Die simultanen {S}ysteme bin\"arer {F}ormen, \href{http://dx.doi.org/10.1007/BF01444021}{\textit{Math. Ann.}}
  \textbf{2} (1870), 227--280.

\bibitem{Gund}
Gundelf\/inger S., Zur Theorie des simultanen Systems einer cubischen und einer
  biquadratischen bin\"aren Form, Habilitationsschrift, J.B.~Metzler,
  Stuttgart, 1869.

\bibitem{Ha}
Hammond J., Note on an {e}xceptional {c}ase in which the {f}undamental
  {p}ostulate of {p}rofessor {S}ylvester's {t}heory of {t}amisage fails,
  \href{http://dx.doi.org/10.1112/plms/s1-14.1.85}{\textit{Proc. London Math. Soc.}} \textbf{14} (1882), 85--88.

\bibitem{Hi2}
Hilbert D., Ueber die vollen {I}nvariantensysteme, \href{http://dx.doi.org/10.1007/BF01444162}{\textit{Math. Ann.}}
  \textbf{42} (1893), 313--373.

\bibitem{HoRo}
Hochster M., Roberts J.L., Rings of invariants of reductive groups acting on
  regular rings are {C}ohen--{M}acaulay, \href{http://dx.doi.org/10.1016/0001-8708(74)90067-X}{\textit{Adv. Math.}} \textbf{13}
  (1974), 115--175.

\bibitem{Mo}
Morley R.K., On the {f}undamental {p}ostulate of {t}amisage, \href{http://dx.doi.org/10.2307/2370111}{\textit{Amer.~J.
  Math.}} \textbf{34} (1912), 47--68.

\bibitem{Olver}
Olver P.J., Classical invariant theory, \href{http://dx.doi.org/10.1017/CBO9780511623660}{\textit{London Mathematical Society
  Student Texts}}, Vol.~44, Cambridge University Press, Cambridge, 1999.

\bibitem{Pro07}
Procesi C., Lie groups. An approach through invariants and representations,
  Universitext, Springer, New York, 2007.

\bibitem{Springer77}
Springer T.A., Invariant theory, \textit{Lecture Notes in Mathematics}, Vol.
  585, Springer-Verlag, Berlin, 1977.

\bibitem{Springer80}
Springer T.A., On the invariant theory of {${\rm SU}_{2}$}, \textit{Indag.
  Math.} \textbf{42} (1980), 339--345.


\bibitem{Sylv3}
Sylvester J.J., D\'etermination du nombre exact des covariants irr\'eductibles
  du syst\`eme cubo-biquadratique binaire, \textit{C.~R.~Acad. Sci. Paris}
  \textbf{87} (1878), 477--481.

\bibitem{Sylv4}
Sylvester J.J., Proof of the hitherto undemonstrated fundamental theorem of
  invariants, \textit{Phil. Mag.} \textbf{5} (1879), 178--188.

\bibitem{Sylv0}
Sylvester J.J., Sur le vrai nombre des covariants \'el\'ementaires d'un
  syst\`eme de deux formes biquadratiques binaires, \textit{C.~R.~Acad. Sci.
  Paris} \textbf{84} (1877), 1285--1289.

\bibitem{Sylv2}
Sylvester J.J., Sur le vrai nombre des formes irr\'eductibles du syst\`eme
  cubo-biquadratique, \textit{C.~R.~Acad. Sci. Paris} \textbf{87} (1878),
  445--448.

\bibitem{Sylv1}
Sylvester J.J., Sur les covariants fundamentaux d'un syst\`eme cubo-quartique
  binaire, \textit{C.~R.~Acad. Sci. Paris} \textbf{87} (1878), 287--289.

\bibitem{SylvF}
Sylvester J.J., Franklin F., Tables of the {g}enerating {f}unctions and
  {g}roundforms for the {b}inary {q}uantics of the {f}irst {t}en {o}rders,
  \href{http://dx.doi.org/10.2307/2369240}{\textit{Amer.~J. Math.}} \textbf{2} (1879), 223--251.

\end{thebibliography}
\end{document}